\documentclass[12pt]{article}
\usepackage{amssymb,amsmath,amsfonts,bbm,pifont,upgreek,bbold,accents}  
\author{Fernando Argentieri }

\font\teneufm=eufm10
\font\seveneufm=eufm7
\font\fiveeufm=eufm5
\newfam\eufmfam
\textfont\eufmfam=\teneufm
\scriptfont\eufmfam=\seveneufm
\scriptscriptfont\eufmfam=\fiveeufm

\newcommand\beq[1]{ \begin{equation}\label{#1} }
\newcommand{\eeq}{ \end{equation} }

\newcommand\beqa[1]{ \begin{eqnarray} \label{#1}}
\newcommand{\eeqa}{ \end{eqnarray} }
\newcommand{\beqano}{ \begin{eqnarray*} }
\newcommand{\eeqano}{ \end{eqnarray*} }

\newtheorem{theorem}{Theorem}
\newtheorem{definition}{Definition}
\newtheorem{proposition}{Proposition}
\newtheorem{lemma}{Lemma}
\newtheorem{sublemma}{Sublemma}
\newtheorem{remark}{Remark}
\newtheorem{notationalremark}{Notations}
\newtheorem{corollary}{Corollary}
\newtheorem{assumption}{Assumption}
\newtheorem{claim}{Claim}
\newtheorem{tools}{$\negsp\negsp$}[subsection]

\newcommand\thm[1]{ \begin{theorem}\label{#1}}
\newcommand\thmtwo[2]{ \begin{theorem}[#1]\label{#2}}
\newcommand\ethm{ \end{theorem} }
\newcommand\dfn[1]{ \begin{definition}\label{#1} \rm}
\newcommand\dfntwo[2]{ \begin{definition}[#1]\label{#2} \rm}
\newcommand\edfn{ \end{definition} }
\newcommand\pro[1]{ \begin{proposition}\label{#1}}
\newcommand\protwo[2]{ \begin{proposition}[#1]\label{#2}}
\newcommand\epro{ \end{proposition} }
\newcommand\lem[1]{ \begin{lemma}\label{#1}}
\newcommand\lemtwo[2]{ \begin{lemma}[#1]\label{#2}}
\newcommand\elem{ \end{lemma} }
\newcommand\sublem[1]{ \begin{sublemma}\label{#1}}
\newcommand\sublemtwo[2]{ \begin{sublemma}[#1]\label{#2}}
\newcommand\esublem{ \end{sublemma} }
\newcommand\rem[1]{ \begin{remark}\label{#1} \rm}
\newcommand\erem{ \end{remark} }
\newcommand\notrem[1]{ \begin{notationalremark}\label{#1} \rm}
\newcommand\enotrem{ \end{notationalremark} }
\newcommand\cor[1]{ \begin{corollary}\label{#1}}
\newcommand\cortwo[2]{ \begin{corollary}[#1]\label{#2}}
\newcommand\ecor{ \end{corollary} }
\newcommand\asmp[1]{ \begin{assumption}\label{#1}}
\newcommand\asmptwo[2]{ \begin{assumption}[#1]\label{#2}}
\newcommand\easmp{ \end{assumption} }
\newcommand\clm[1]{ \begin{claim}\label{#1}}
\newcommand\eclm{ \end{claim} }
\newcommand{\proof}{\par\medskip\noindent{\bf Proof\ }}
\newcommand\equ[1]{{\rm (\ref{#1})}}
\expandafter\chardef\csname pre amssym.def
at\endcsname=\the\catcode`\@
\catcode`\@=11
\def\undefine#1{\let#1\undefined}
\def\newsymbol#1#2#3#4#5{\let\next@\relax
 \ifnum#2=\@ne\let\next@\msafam@\else
 \ifnum#2=\tw@\let\next@\msbfam@\fi\fi
 \mathchardef#1="#3\next@#4#5}
\def\mathhexbox@#1#2#3{\relax
 \ifmmode\mathpalette{}{\m@th\mathchar"#1#2#3}%
 \else\leavevmode\hbox{$\m@th\mathchar"#1#2#3$}\fi}
\def\hexnumber@#1{\ifcase#1 0\or 1\or 2\or 3\or 4\or 5\or 6\or 7\or
8\or
 9\or A\or B\or C\or D\or E\or F\fi}
\ifcase\@ptsize
 \font\tenmsb=msbm10
 \font\sevenmsb=msbm7
 \font\fivemsb=msbm5
\or
 \font\tenmsb=msbm10 scaled \magstephalf
 \font\sevenmsb=msbm7 scaled \magstephalf
 \font\fivemsb=msbm5  scaled \magstephalf
\or
 \font\tenmsb=msbm10 scaled \magstep1
 \font\sevenmsb=msbm7 scaled \magstep1
 \font\fivemsb=msbm5 scaled \magstep1
\fi
\newfam\msbfam
\textfont\msbfam=\tenmsb
\scriptfont\msbfam=\sevenmsb
\scriptscriptfont\msbfam=\fivemsb
\edef\msbfam@{\hexnumber@\msbfam}
\def\Bbb#1{\fam\msbfam\relax#1}
\def\widehat#1{\setboxz@h{$\m@th#1$}%
 \ifdim\wdz@>\tw@ em\mathaccent"0\msbfam@5B{#1}%
 \else\mathaccent"0362{#1}\fi}
\def\widetilde#1{\setboxz@h{$\m@th#1$}%
 \ifdim\wdz@>\tw@ em\mathaccent"0\msbfam@5D{#1}%
 \else\mathaccent"0365{#1}\fi}

\def\RIfM@{\relax\ifmmode}
\def\nonmatherr@#1{\errmessage{\string#1\space allowed only in math mode}}
\def\Bbb{\RIfM@\expandafter\Bbb@\else
 \expandafter\nonmatherr@\expandafter\Bbb\fi}
\def\Bbb@#1{{\Bbb@@{#1}}}
\def\Bbb@@#1{\fam\msbfam\relax#1}
\def\setboxz@h{\setbox\z@\hbox}
\def\wdz@{\wd\z@}
\catcode`\@=\csname pre amssym.def at\endcsname
\newcommand{\giu}{{\medskip\noindent}}
\newcommand{\Giu}{{\bigskip\noindent}}
\newcommand{\nl}{{\smallskip\noindent}}

\newcommand{\qed}{\hskip.5truecm
\vrule width 1.7truemm height 3.5truemm depth 0.truemm
\par\Giu}

\newcommand{\negsp}{\hspace{-.09truecm}}  

\newcommand{\IIt}{\mathcal{I}_{\t}}
\newcommand{\DDt}{\mathcal{D}_{\t}}

\newcommand{\pq}{\frac{p}{q}}

\newcommand{\Dgt}{D_{\gamma,\tau}}
\newcommand{\gta}{\gamma(\alpha,\tau)}
\newcommand{\gnat}{\g_{n}(\a,\t)}
\newcommand{\gtal}{\g_{-}(\a,\t)}
\newcommand{\gttal}{\g_{-}(\a,\t')}
\newcommand{\gtar}{\g_{+}(\a,\t)}
\newcommand{\gttar}{\g_{+}(\a,\t')}
\newcommand{\gmat}{\g_{m}(\a,\t)}

\newcommand\ugper[1]{ \stackrel{#1}{=} }

\newcommand{\dst}{\displaystyle}

\newcommand\su[1]{ \frac{1}{ {#1}} }

\renewcommand{\natural}{ {\Bbb N}   }
\newcommand{\real}{ {\Bbb R}   }

\renewcommand{\a }{ {\alpha}   }
\renewcommand{\b}{ {\beta}   }
\newcommand{\g}{ {\gamma}   }

\renewcommand{\d}{ {\delta}   }
\newcommand{\D}{ {\Delta}   }
\newcommand{\e }{ {\epsilon}   }

\newcommand{\m}{ {\mu}   }

\renewcommand{\t}{ {\tau}   }

\renewcommand\subset{\subseteq}

\begin{document}
\title{Isolated points of Diophantine sets}
\maketitle
\begin{abstract}
  Let $\g\in(0;\su{2}),\t\geq 1$ and define the ``$\g,\t$ Diophantine set" as: 
  $$\Dgt:=\{\a\in (0;1): ||q\a||\geq\frac{\g}{q^{\t}}\quad\forall q\in\Bbb{N}\},\qquad ||x||:=\inf_{p\in\Bbb{Z}}|x-p| $$
  We analyze the topology of these sets and we show that generally they have isolated points.
\end{abstract}
\section{Introduction}
Diophantine sets play an important role in dynamical systems, in particular, in
small divisors problems with applications to KAM theory, Aubry-Mather theory,
conjugation of circle diffeomorphisms, etc. (see, for example, [3], [5], [9], [12], [13], [14], [16]). 

\nl
The set $\Dgt$ is compact and totally
disconnected (since $\Dgt\cap\Bbb{Q}=\emptyset$), however,  it is not clear whether, for some $\gamma$ and $\tau$, there exist isolated points in  $\Dgt$.

\giu
In this paper, we provide explicit examples of $\Dgt$ with isolated points, giving, in particular, a partial answer to a 
a question raised by Broer in [2] (see remark (iii) below). 

\nl
Our main results are the following.

\giu
{\bf {Proposition 1}}
{\sl 
Let $n\in \natural$, $n\geq 2$ and define 
\beq{gt}
\bar{\a}:=\frac{\sqrt{n^2+4}-n}2\ ,\ \qquad \g:=\frac{1}{\bar{\a}+n}\ ,\qquad\ \t:=\frac{\log (\bar{\a}+n)}{\log n}\ .
\eeq
Then $\bar{\a}$ is an isolated point of $D_{\g,\t}$.}

\Giu

\nl
Indeed, we can show that, for all Diophantine numbers, there exists an `equivalent number' that is isolated in some Diophantine set:

\nl
{\bf{Theorem A}}
{\sl 
Let $\g\in(0,\su{2})$, $\t\geq 1$. Define the map:
\beq{}
\Phi_{\g,\t}(z):=\frac{\eta z+1}{(2\eta+1)z+2}
\eeq
with
\beq{}
\eta:=\Big[\frac{2^{\t} 3}{\g}\Big].
\eeq
Then $\Phi(\Dgt)\subset D_\tau:=\bigcup_{\gamma>0} \Dgt$. Moreover, for all $\a\in\Phi(\Dgt)$ there exists $\t_{\a}>\t$ and $\gamma_\a>0$ such that $\a$ is isolated in  $D_{\gamma_\alpha, \tau_\alpha}$.
}

\giu
The isolated points constructed in Theorem A  depend only on the first coefficients of their continued fraction (that we can change up to an equivalent number). 

\giu
Finally,  we show that a Diophantine number may be an isolated point ``for infinitely many $\t$":

\giu
{\bf{Theorem B}} {\sl
Fix $\t\geq 1$ and a strictly decreasing sequence $\{\t_{n}\}_{n\in\Bbb{N}}$ with
$\t_n\searrow \t$. Then, there exist $\g>0$, $\a\in \Dgt$ and  sequences $\{{\bar{\t}}_{n}\}_{n\in\Bbb{N}}$, $\{\g_n\}_{n\in\Bbb{N}}$ with, $\bar{\t}_n\in(\t_n,\t_{n+1})$, $\g_n\searrow \g$ such that $\a$ is an isolated point of $D_{\g_n,\bar{\t}_{n}}$ for all $n$.
}

\giu
{\bf Remarks} (i) The existence of isolated points of Diophantine sets may be related to isolated tori and KAM stability
in two degrees of freedom.

\nl
(ii) Our analysis is based on continued fractions and relations with dynamics in higher dimensions are, therefore, not
clear.

\Giu
(iii) The paper [2] is entitled: ``{\sl Do Diophantine vectors form a Cantor bouquet?}'', namely, is the set
$\D^{N}_{\g,\t}\cap \Bbb{S}^{N-1}$, where 
$$\D^{N}_{\g,\t}:=\{\omega\in\Bbb{R}^{N}:|\omega\cdot n|\geq\frac{\g}{|n|^{\t}}\quad\forall n\in\Bbb{Z}^{N},n\not=0\}\,,$$
and $\Bbb{S}^{N-1}$ denotes the unit sphere in $\real^N$, a Cantor set?

\nl
In dimension $N=2$ it is clearly equivalent  to consider the intersection of $\D^{2}_{\g,\t}$ with the line $\omega_2=1$, which, upon restricting to the unit interval, coincides with the set $\Dgt$.

\nl
{\sl Our results, therefore,  show that, in general, the answer to such a question is negative, at least, in dimension $N=2$}.

\Giu
(iv) In all our examples of isolated points the following holds: if $\a$ is an isolated point of $\Dgt$, then $\g$ is the best constant such that the Diophantine conditions with exponent $\t$ holds. By an amazing Theorem of Roth, for any algebraic numbers $\a$, given $\t>1$ there exists $\g>0$ such that $\a\in\Dgt$ (see, for example,  [11]).
We believe that, for algebraic numbers of degree greater then 2,  the statement of Theorem B holds. So, information about isolated points may be in connection with continued fraction properties of algebraic numbers.

\Giu
The paper is organized as follows: in section 2 we give the main definitions and make few remarks, in section 3 we remind general properties of Diophantine sets, in section 3 we  prove our main results. Finally, section 5 contains a few concluding  observations and some questions.

\section{Definitions and remarks}
\nl
\subsection{Definitions}

\begin{itemize}
    \item $\Bbb{N}:=\{1,2,3,...\}$, $\Bbb{N}_{0}:=\{0,1,2,3,...\}$ 
    \item Given $a,b\in\Bbb{Z} -\{0\}$, we indicate with $(a,b)$ the maximum common divisor of $a$ and $b$.
    \item Let $\a$ be a real number. We indicate with $[\a]$ the integral part of $\a$, with $\{\a\}$ the fractional part of $\a$ .
     \item Given E$\subset{\Bbb{R}}$, we indicate with $\mathcal{I}$(E) the set of isolated points of E.
     \item Given E$\subset{\Bbb{R}}$, we indicate with $\mathcal{A}$(E) the set of accumulated points of E.
     \item We say that E$\subset{\Bbb{R}}$ is perfect if $\mathcal{A}$(E)=E.
    \item Given a Borel set E$\subset{\Bbb{R}}$ we denote with $\m$(E) the Lebesgue measure of E.
    \item A topological space X is a totally disconnected space if the points are the only connected subsets of X.
    \item $X\subset\Bbb{R}$ is a Cantor set if it is closed, totally disconnected and perfect.
    \item For $E\subset\Bbb{R}^{n}$, $\dim_{H}E$ is the Hausdorff dimension of $E$.
    \item Given $\a\in\Bbb{R}$ we define:
    $$||\a||:=\min_{p\in\Bbb{Z}}|\a -p|$$
    \item Given $\g>0, \t\geq1$, we define the $(\g,\t)$ Diophantine points in $(0;1)$ as the numbers in the set:
    $$ \Dgt:=\{\a\in (0 ;1): ||q\a||\geq\frac{\g}{q^{\t}}\quad\forall q\in\Bbb{N}\}$$
    \item $$ D^{\Bbb{R}}_{\g,\t}:=\{\a\in\Bbb{R}:||q\a||\geq\frac{\g}{q^{\tau}}\quad\forall q\in\Bbb{N}\},$$ $$D_{\t}:=\bigcup_{\g>0}D_{\g,\t},\quad D:=\bigcup_{\t\geq 1}D_{\t}.$$  We call $D$ the set of Diophantine numbers.
    \item Given $\t\geq 1,\a\in \Bbb{R}$, we define:
    $$\gta:=\inf_{q\in\Bbb{N}}q^{\t}||q\a||$$
    \item Given $\a\in\Bbb{R}$ we define:
    $$\t(\a):=\inf\{\t\geq 1:\gta>0\}$$
   \item Given an irrational number $\a=[a_{0};a_{1},...]:=a_0+\su{a_1+\su{a_2+...}}$, we denote with $\{\frac{p_n}{q_n}\}_{n\in\Bbb{N}_{0}}$ the convergents of $\a$, $\a_{n}:=[a_{n};a_{n+1},...]$\footnote{for information about continued fractions see [4],[8],[15] }. 
   \item We indicate with $[a_1,a_2,a_3,...]:=\su{a_1+\su{a_2+\su{a_3+...}}}$.
   \item Let $\a$ be an irrational number. We define:
   $$\gnat:=q_n^{\t}||q_n\a||=q_n^{\t}|q_n\a-p_n|$$
   \item Let $\t\geq 1$,
   $$\gtal:=\inf_{n\in2 \Bbb{N}_{0}}\gnat,$$
   $$\gtar:=\inf_{n\in2\Bbb{N}_{0}+1}\gnat,$$
   $${\DDt}:=\{\a\in D_{\t} :\t(\a)=\t\},$$
   $${\IIt}:=\{\a\in D_{\t}: \exists n\not\equiv m\quad{(\rm{mod} 2)}, \gnat=\gmat=\g(\a,\t)\}.$$
   $${\mathcal{I}}:=\cup_{\t\geq 1}{\IIt}$$
   \item Let $p\in\Bbb{Z},q\in\Bbb{N}$, $\g>0,\t\geq 1$. We define: $I_{\g,\t}(p,q):=\left(\frac{p}{q}-\frac{\g}{q^{\t+1}};\pq+\frac{\g}{q^{\t+1}}\right)$.

\end{itemize}
\subsection{Remarks}

\begin{enumerate}
    \item $\a\in \Dgt\iff 1-\a\in\Dgt$.
    \item $\gta\leq \min\{\a,1-\a\}.$
    \item Fixed $\t\geq 1$, $\g(.,\t):D_{\t} \rightarrow (0,\frac{1}{2})$.
    \item $\Dgt^{\Bbb{R}}=\bigcup_{n\in\Bbb{Z}}(\Dgt+n)$, thus we can restrict to study the Diophantine points in $(0,1)$.
    \item 
    \beq{fond}
    \left\{
    \begin{array}{l}
           \gnat=\frac{q_{n}^{\t}}{\a_{n+1}q_{n}+q_{n-1}},\\
           \su{\gnat}=\frac{q_{n+1}}{q_{n}^{\t}}+\frac{1}{\a_{n+2}q_{n}^{\t-1}}
    \end{array}\right.
    \eeq
    \item $\gta=\inf_{n\in\Bbb{N}_{0}}\gnat$.
    \item If $\t<\t(\a)$, then $\gta=0$; if $\t>\t(\a)$ then $\gta>0$. Moreover, for $\t>\t(\a)$ the inf is a minimum.
    \item $\a\in {\DDt}\iff \t(\a)=\t$ and $\gta>0$.
    \item If $\a\in{\IIt}$, then $\a$ is an isolated point of $\Dgt$.
    \item The cardinality of ${\IIt}$ is at most countable.
    \item $\m({\DDt})=0$ for all $\t\geq 1$.
    \item $\g_{0}(\a,\t)=\left\{\a\right\}$, in particular $\g_{0}(\a,\t)$ does not depend on $\t$.
    \item Let $\pq$ a rational number. 
    \beq{}
    \a\in D_{\t}\iff \left\{\a+\pq\right\}\in D_\t,
    \eeq
    \beq{}
    \a\in{\DDt}\iff \left\{\a+\pq\right\}\in{\DDt}.
    \eeq
    \item If $\t>\t(\a)$, $\gtal=\gtar$, then $\a\in{\IIt}$.
    \item $\a\in D_\t\iff q_{n+1}=O(q_{n}^{\t}).$
\end{enumerate}
\proof
(a), (d) are clear, (b) follows by definition of $\g(\a,\t)$ and by remark (a). (c) follows by (b) ($\a$ is in $(0;1)$).

\giu
(e): the first formula follows by properties of continued fractions, moreover:
\beq{}
\su{\gnat}=\frac{\a_{n+1}q_{n}+q_{n-1}}{q_{n}^{\t}}=\frac{(a_{n+1}q_{n}+q_{n-1})+\frac{q_n}{\a_{n+2}}}{q_{n}^{\t}}\\
=\frac{q_{n+1}}{q_{n}^{\t}}+\frac{1}{\a_{n+2}q_{n}^{\t-1}}.
\eeq

\giu
(f): follows by:
\beq{}
||q_n\a||=\min_{1\leq q\leq q_n}||q\a||
\eeq
and by definition of $\g(\a,\t)$.

\giu
(g): The first part is clear. To prove that for $\t>\t(\a)$ the inf is a minimum, take $\t(\a)<\t'<\t$, then $\g(\a,\t')>0$ and 
\beq{infinito}
\lim_{n\rightarrow+\infty}q_n^{\t}||q_n\a||=\lim_{n\rightarrow+\infty}q_n^{\t-\t'}q_n^{\t'}||q_n\a||\geq\lim_{n\rightarrow+\infty}\g(\a,\t')q_n^{\t-\t'}=+\infty.
\eeq
By (\ref{infinito}) there exists $N\in\Bbb{N}$ such that for all $n\geq N$, 
\beq{}
\gnat>\g(\a,\t).
\eeq
Therefore the inf is reached and it is a minimum.

\giu
(h): It is obvious.

\giu
(i): If $\a$ is in ${\IIt}$, there exist $n$ even and $m$ odd such that:
\beq{}
\g(\a,\t)=\gnat=\g_m(\a,\t).
\eeq
So $\a$ is separated by the two intervals $I_{\g,\t}\left(p_n,q_n\right)$ and $I_{\g,\t}\left(p_m,q_m\right)$. Then, noting that $I_{\g,\t}\left(p,q\right)\subset\Dgt^{c}$ for all $p\in\Bbb{Z},q\in\Bbb{N}$, we get (i).

\giu
(j): If $\gnat=\gmat=\gta$, with $n$ even, $m$ odd, then:
\beq{num}
\a=\frac{p_n}{q_n}+\frac{\frac{p_m}{q_m}-\frac{p_n}{q_n}}{1+(\frac{q_n}{q_m})^{\t+1}},
\eeq
\beq{}
\g=\frac{\frac{p_m}{q_m}-\frac{p_n}{q_n}}{\su{q_{n}^{\t+1}}+\su{q_{m}^{\t+1}}},\eeq
so ${\IIt}$ is at most countable.

\giu
(k): $\m(D_1)=0$ ($D_1$ is the set of numbers with bounded coefficients of the continued fraction). Moreover $\m(D_\t)=1$ for all $\t>1$ (because of $\m(\Dgt^c)=O(\g)$). For $1<\t'<\t$ we have ${\mathcal{D}_\t}\subset {D_{\t'}^{c}}$. So, for $\t>1$: $\m\left({\mathcal{D}_\t}\right)=0$.

\giu
(l), (m): They are obvious.

\giu
(n): Because of $\t>\t(\a)$, as in the proof of (g) we get that $\gtal$ and $\gtar$ are reached, so there exist $n$ even and $m$ odd with $\gtal=\gnat$, $\gtar=\g_m(\a,\t)$. Now (n) follows by definition of ${\IIt}$ and by 
$\gtal=\gtar=\gta$.

\giu
(o): It follows by (e) and (f).
\section{Basic properties of Diophantine sets}
Let us recall some simple facts about Diophantine sets. The case $\t=1$ is quite different to the others. 
\rem{}
If $0<\g'\leq \g$, $\t'\geq\t\geq1$, then $\Dgt\subset D_{\g',\t'}$. Moreover, $\Dgt$ is compact and totally disconnected (because of $\Dgt\cap\Bbb{Q}=\emptyset$).
\erem{}
\rem{}
$D_1$ is the set of irrational numbers with bounded coefficients of their continued fractions.
\erem{}
\proof
It follows by (\ref{fond}). \qed
\thm{}{\bf (Hurwitz)}(see \cite{8})
Let $\a$ be an irrational number. There exist infitely many $q\in\Bbb{N}$ such that 
\beq{}
q||q\a||<\su{\sqrt{5}q}.
\eeq
\ethm{}
\thm{}{\bf (Borel)}(see \cite{1})
\label{t6}
Given a function $\psi:\Bbb{N}\rightarrow\Bbb{N}$, define
$$A(\psi):=\{[a_{0};a_{1},...,a_{n},...]:0<a_{n}<\psi(n)\}.$$
Then:
\beq{}
\sum_{n\in\Bbb{N}}\frac{1}{\psi(n)}<\infty\Rightarrow\quad \m(A)>0,
\eeq
\beq{}
\sum_{n\in\Bbb{N}}\frac{1}{\psi(n)}=\infty\Rightarrow\quad \m(A)=0.
\eeq
\ethm{}
\rem{}
By Hurwitz's theorem, if $\g>\su{\sqrt{5}}$, then $D_{\g,1}=\emptyset$.
\erem{}
\rem{}
For all $\g\in(0,\su{2})$ we have $\m(D_{\g,1})=0$. In particular $\m(D_1)=0$.
\erem{}
\proof
It follows by (\ref{fond}) and Borel's theorem. \qed
\giu

\giu
Unless $D_1$ has zero measure, it has positive Hasdorff dimension. In fact, the following holds:
\thm{} {\bf (Jarnik)}(see \cite{10})
$\dim_{H}(D_1)=1$.
\ethm{}
\thm{}(see \cite{6})
Let $\g>\su{3}$. Then the set:
\beq{}
\left\{\a\in(0,1):\liminf q||q\a||\geq \g\right\}
\eeq
is at most countable. In particular, for $\g>\su{3}$ $D_{\g,1}$ is at most countable.
\ethm{}
The case $\t>1$ is quite different.
\rem{}
Let $\t>1$. Then, for $\g>0$ we have 
\beq{}
\m(\Dgt^{c})=O(\g).
\eeq
In particular, $\m(D_\t)=1$ for all $\t>1$.
\erem{}
\proof
For $\t>1$:
\beq{}
\m(\Dgt^{c})\leq\sum_{q\in\Bbb{N}}\sum_{0\leq p\leq q-1}\frac{2\g}{q^{\t+1}}=2\g\sum_{q\in\Bbb{N}}\su{q^{\t}}=O(\g).
\eeq
\qed
\cor{}
\beq{}
\m\left(\bigcap_{\t>1}D_{\t}\right)=1.
\eeq
\ecor{}

\section{Isolated points of Diophantine sets}
In this section we give the proof of the results. We start by proving Proposition 1.
\proof
\Giu 
Fix $\a:=\bar{\a}+n$. It is easy to verify that $\a$ is such that:
\beq{alpha}
\left\{
\begin{array}{l}
\dst \a=\su{\a}+n\ , \qquad n^\t=\a\ ,\\  \ \\
\dst \a=[n;n,n,n,....]:=n+\su{n+\su{n+...}}\ ,\\ \ \\
p_0=n,\ q_0=1,\ p_1=n^{2}+1,\  q_1=n,\  \a_k=\a\quad\forall k\geq 1,\ q_{k+1}=p_k\ (\forall k\ge 0)\ .
\end{array}\right.
\eeq
\Giu

\Giu

\Giu
\nl
For $k=0$:
\beq{a0}
\Big| \a- \frac{p_0}{q_0}\Big|\ugper{\equ{alpha}} \a-n\ugper{\equ{alpha}} \su{\a}\ugper{\equ{gt}}\g .
\eeq
For $k\ge 1$, from \equ{alpha} and the fact that $p_k/q_k\le p_1/q_1$ and $q_k\ge q_1$, we obtain:
\beqano
\frac{q_{k+1}}{q_{k}^{\t}}+\frac{1}{a_{k+2}q_{k}^{\t-1}}&\ugper{}&\frac{p_{k}}{q_{k}}\frac{1}{q_{k}^{\t-1}}+\frac{1}{\a q_{k}^{\t-1}}\\
&\leq& \frac{p_{1}}{q_{1}}\frac{1}{q_{1}^{\t-1}}+\frac{1}{\a q_{1}^{\t-1}}\ugper{}\frac{n^2+1}{n^\t}+\su{n^{\t-1}\a}\\
&\ugper{}& \frac{n^2+1}{\a} +\frac{n}{\a^2}=\su\a \Big(n^2+1+\frac{n}{\a}\Big)\\
&=&\su\a \ (\a n+1)=n+\su\a =\a\\
&=&\su\g\ ,
\eeqano
that, togheter with \equ{a0},  it shows that $\a\in D_{\g,\t}+n$. 
\\
From \equ{alpha}, 
\beqano
\Big| \a -\frac{p_1}{q_1}\Big|&=& \frac{p_1}{q_1} - \a= \frac{n^2+1}{n}-\a=\su{n}+n-\a\\
&=&\su{n}-\su\a=\su{n\a^2}=
\su\a\ \su{q_1 n^\t}=\su{\a q_1^{\t+1}}\\
&=&\frac{\g}{q_1^{\t+1}}\ ,
\eeqano
that shows, togheter with \equ{a0}, that $\a$ divides the two intervals  $I_{\g,\t}(p_0,q_0)$ and $I_{\g,\t}(p_1,q_1)$, with $I_{\g,\t}(p,q):=\left(\pq-\frac{\g}{q^{\t+1}};\pq+\frac{\g}{q^{\t+1}}\right)$. 
So $\a\in D_{\g,\t}+n$ implies that $\a$ is an isolated point of $D_{\g,\t}+n$, i.e. $\bar{\a}$ is an isolated point of $\Dgt$. \qed

\nl
Before proving Theorem A we need some simple lemma. So we prove at first the continuity of the functions $\gta$, $\gtal,\gtar$ as functions of $\t$.

\lem{gen}
Let $a\in\Bbb{R}$, $f_n\geq 0$ be continuous and increasing functions in $[a,+\infty)$ such that:
\nl
\beq{aaa}
 \forall x>a,\quad \lim_{n\rightarrow+\infty} f_n(x)=+\infty.
\eeq
Define
\beq{ccc}
f(x):=\inf_{n\in\Bbb{N}}f_{n}(x).
\eeq
If $f$ is bounded,  then $f\in C([a,+\infty))$ and $f$ is an increasing function.

\elem{}
\proof
Observe that $f$ is increasing because $f_n$ are increasing. Let $C>0$ be such that $f(x)\leq C$ for all $x\in[a,+\infty)$.
Take $x\in\Bbb{R}$ such that $a<x$. By (\ref{aaa}) there exists $N\in\Bbb{N}$ such that for all $n\geq N$, $f_{n}(x)>C>0$. For $y\geq x$, $f(y)=\min_{0\leq n<N} f_{n}(y)$, so $f$ is continuous and increasing in $(x,+\infty)$ and $f\in C((a,+\infty))$. It remains to show that $f$ is continuous in $a$, i.e. $f(a)=\lim_{x\rightarrow a}f(x)$. In fact, for all $\epsilon>0$ there exists $n\in\Bbb{N}$ such that 
\beq{}
0<f_{n}(a)-f(a)<\epsilon
\eeq
and by continuity of $f_n$ there exists $\delta>0$ such that for $0<x-a<\delta$ we have:
\beq{}
0<f_{n}(x)-f_{n}(a)<\epsilon.
\eeq
So, for $0<x-a<\delta$:
\beq{}
0\leq f(x)-f(a)\leq f_{n}(x)-f_{n}(a)+f_{n}(a)-f(a)<2\e,
\eeq
that proves the continuity in $a$.
\cor{a}
Fixed $\a\in D$, the functions $\g(\a,\t),\gtal,\gtar$ are continuous and increasing for $\t\geq \t(\a)$.
\ecor{}
\proof
We prove the corollary for $\gta$ (the proof for $\gtal,\gtar$ are similar). Observe that $\gnat\leq \su{2}$. 
Consider the $\gnat$ as functions of $\t$. For $\t>\t(\a)$ we have
\beq{}
\lim_{n\rightarrow +\infty}\gnat=+\infty
\eeq
Moreover the $\gnat$ are increasing with respect to $\t$, so the hypothesis of Lemma \ref{gen} are satisfied. \qed
Now we give a simple sufficient condition such that a Diophantine number belongs to ${\IIt}$ for some $\t\geq\t(\a)$.

\lem{gen}
Let $\a\in D\cap(0;\su{2})$ be such that there exists $\t'>\t(\a)$ with:
\beq{cond}
\gttal\geq\gttar
\eeq
Then there exists $\t\geq\t'$ such that $\a\in{\IIt}$
\elem{}
\proof
If:
\beq{}
\gttal=\gttar
\eeq
then $\a\in{\mathcal{I}}_{\t'}$ by remark (g) and because of $\t'>\t(\a)$. Now consider the case:
\beq{}
\gttal>\gttar
\eeq
Observe that:
\beq{m}
\gtal\leq \g_{0}(\a,\t)\leq\max\{\a,1-\a\}.
\eeq
Moreover
\beq{mm}
\lim_{\t\rightarrow+\infty}\gtar=+\infty
\eeq
because it is an increasing function and because of $\a\in(0,\su{2})$.
So, by continuity of $\gtal,\gtar$ and by (\ref{m}), (\ref{mm}) we get that there exists $\t>\t'$ such that $\gtal=\gtar$, so $\a\in{\IIt}$ by remark (g).\qed
\rem{}
Note that the condition (\ref{cond}) is satisfied for $\bar{\a}$ defined in the Proposition.
Moreover, for this $\bar{\a}$ there exists a unique $\t$ such that $\g_{-}(\b,\t)=\g_{+}(\b,\t)$.
\erem{}
{\bf{Proof (Theorem A)}}
Fixed $\t\geq1,\g\in(0;\su{2})$, consider the map $\Phi_{\g,\t}$ defined in the statement of Teorem A. Let $\a\in\Dgt$. Observe that, if $\a=[a_1,a_2,...]$ then:
\beq{}
\Phi(\a)=[2,\left[2^{\t}\frac{3}{\g}\right],a_1,a_2,...].
\eeq
We denote with $q_n$ the denominator of the n-th convergent to $\Phi(\a)$, with $\b_n$ the n-th residue of $\Phi(\a)$ and with $q_n'$ the denominator of the n-th convergent to $\a$.
We recall that:
\beq{}
\su{\g_n(\Phi(\a),\t)}=\frac{q_{n+1}}{q_{n}^{\t}}+\frac{1}{\b_{n+2}q_{n}^{\t+1}},
\eeq
and 
\beq{z}
\frac{q_{n+1}}{q_{n}^{\t}}+\frac{1}{\b_{n+2}q_{n}^{\t+1}}=\frac{q_{n-1}}{q_{n}^{\t}}+\frac{a_{n+1}}{q_{n}^{\t-1}}+\su{\b_{n+2}q_{n}^{\t+1}}.
\eeq
So, by (\ref{z}):
\beq{beta}
\left\{
\begin{array}{cc}
\dst \su{\g_0(\Phi(\a),\t)}<\su{\Big[\frac{\g}{3}\Big]} \\
\dst \su{\g_1(\Phi(\a),\t)}>\su{\Big[\frac{\g}{3}\Big]}\\
\dst \su{\g_n(\Phi(\a),\t)}<\frac{2}{\g}\quad for\quad n\geq 2.\\ \\
\end{array}\right.
\eeq
In fact:
\beq{}
\su{\g_{0}(\Phi(\a),\t)}=q_1+\su{\b_2}=2+\su{\b_2}<3<\su{\Big[\frac{\g}{3}\Big]},
\eeq
\beq{}
\su{\g_1(\Phi(\a),\t)}>\frac{q_2}{q_1^{\t}}=\frac{2\Big[2^{\t}\frac{3}{\g}\Big]+1}{2^{\t}}\geq\su{\Big[\frac{\g}{3}\Big]},
\eeq
while, for $n\geq 2$:
\beq{}
\su{\g_n(\Phi(\a),\t)}=\frac{q_{n-1}}{q_{n}^{\t}}+\frac{a_{n-1}}{q_{n}^{\t-1}}+\su{\a_{n-2}q_{n}^{\t+1}}<
\eeq{}
\beq{}
<1+\frac{a_{n-1}}{q_{n-2}'^{(\t-1)}}<1+\su{\g}<\su{\Big[\frac{\g}{2}\Big]},
\eeq
using in the first inequality that $q_n>q_n'>q_{n-2}'$.
By (\ref{beta}), for all $\a\in\Dgt$, $\Phi(\a)$ satisfies the hypothesis of Lemma 2.
In fact the first coefficient of $\Phi(\a)$ is greater then $1$, moreover:
\beq{}
\g_{-}(\Phi(\a),\t)>\Big[\frac{\g}{3}\Big]>\g_{+}(\Phi(\a),\t).
\eeq
\nl
So, given $\a\in\Dgt$, $\Phi(\a)$ is a Diophantine number equivalent to $\a$ that is in ${\mathcal{I}_{\t'}}$ for some $\t'>\t$.
From the arbitrariness of $\g,\t$, Theorem A follows. \qed 
\cor{}
For all $\t\geq 1$ we have:
\beq{}
\m\Big(\bigcup_{\t'\geq\t}{\mathcal{I}_{\t'}}\Big)>0.
\eeq
\ecor{}
\proof
It suffices to note that for all $\g\in(0,\su{2}),\t\geq 1$, the map: $\Phi:\Dgt\rightarrow D$ is Lipschitz and that $\m(\Dgt)>0$ for small $\g$. \qed
\rem{iso}
Suppose that $\a\in D$ such that $\gtal=\gtar$ for some $\t>\t(\a)$. Then $\a$ is an isolated point of $D_{\gta,\t}$. 
\erem{}
\proof
In fact, for $\t>\t(\a)$ $\gtal$ and $\gtar$ are achieved for some $n$ even and $m$ odd. \qed
\rem{strano}
If $\gtal=\gtar$ with $\a\in D$ and $\t=\t(\a)$, in general $\a$ is not an isolated point of $D_{\gta,\t}$.
\erem{}
\proof
For example, take $\t=2,\g=\su{4}$. We define $\a=[a_1,a_2,...]$ iteratively.
$a_1:=2$, and for $n\geq 1$:
\beq{}
a_{n+1}:=\frac{q_{n}^{\t-1}}{\g}-3
\eeq
with $q_{-1}=0,q_{0}=1$, $q_n=a_{n-1}q_{n-1}+q_{n-2}$ for $n\geq1$.
Then it is easy to check that the $a_n$ are strictly increasing, moreover $\t(\a)=\t=2$, $\g(\a,\t(\a))=\g=\su{4}$.
For $n\geq 2$ define:
\beq{}
\d_n:=[a_1,a_2,...,a_{n-1},a_{n}+1,1,1,1,...].
\eeq
We show that $\d_k\in \Dgt$ and $\d_k\rightarrow \a$.
For $n<k-1$ we have:
\beq{}
\su{\g_n(\d_k,\t)}<\frac{a_{n+1}}{q_{n}^{\t-1}}+\frac{q_{n-1}}{q_{n}^{\t}}+\frac{1}{q_{n}^{\t-1}}=\su{\g}+\frac{q_{n-1}}{q_{n}^{\t}}-\frac{2}{q_{n}^{\t-1}}<\su{\g}
\eeq
For $n>k-1$ it is clear that
\beq{}
\su{\g_{n}(\d_k,\t)}<2.
\eeq
For $n=k-1$:
\beq{}
\su{\g_n(\d_k,\t)}<\frac{a_{n+1}+1}{q_{n}^{\t-1}}+\frac{q_{n-1}}{q_{n}^{\t}}+\frac{1}{q_{n}^{\t-1}}=\su{\g}+\frac{q_{n-1}}{q_{n}^{\t}}-\frac{1}{q_{n}^{\t-1}}<\su{\g}
\eeq
So we have proved that $\d_k\in\Dgt$ for all $k\geq 2$. Moreover $\d_k\rightarrow\a$, so $\a$ is not an isolated point of $\Dgt$. \qed
The number constructed in the proof of Remark (\ref{strano}) is not an isolated point because the sequence $\su{\gnat}$ converges too slowly to $\su{\g}$. Moreover, observe that $\g(\a,\t)$ is not achieved  ($\gnat<\g$ for all $n$).  

\proof {\bf{(Theorem B)}}
We construct $\a=[a_1,a_2,...]$ with $a_n$ defined iteratively.
We fix:
\beq{}
\left\{
\begin{array}{l}
\dst a_1=3, \quad a_2= \Big[3^{\t_1+1}\Big],\\ \ \\
\dst q_0=1,\quad q_1=a_1,\quad q_2=a_1 a_2+1\\
\end{array}
\right.
\eeq
Define:
\beq{}
C_1:=\max_{k=0,1}\frac{q_{k+1}}{q_{k}^{\t_2}}=\frac{q_{2}}{q_{1}^{\t_2}}>3.
\eeq
For $n\geq 3$ let:
\beq{}
b_n^{(1)}:=\Big[({C_1}^{2}q_{n-1})^{\t_2-1}\Big].
\eeq
As long as $n$ is even or 
\beq{}
\frac{b_n^{(1)}}{q_{n-1}^{\t_1-1}}\geq C_1-1,
\eeq
define 
\beq{}
a_n=1.
\eeq
Because of $q_{n-1}>2^{n-1}$ and $\t_1>\t_2$, there exists $n_1$ such that:
\beq{}
\frac{b_{n_{1}}^{(1)}}{q_{n_1-1}^{\t_1-1}}<C_1-1.
\eeq
For such $n_1$, define 
\beq{}
a_{n_{1}}=b_{n_{1}}.
\eeq
Define:
\beq{}
C_2:=\max_{k\leq n_1}\frac{a_k}{q_{k-1}^{\t_3-1}}=\frac{a_{n_{1}}}{q_{n_{1}-1}^{\t_3-1}}>{C_1}^{2}-1
\eeq
For $n>n_1$, define:
\beq{}
b_n^{(2)}:=\Big[({C_2}^{2}q_{n-1})^{\t_3-1}].
\eeq
As long as $n$ is odd or
\beq{}
\frac{b_n^{(2)}}{q_{n-1}^{\t_2-1}}\geq C_2-1
\eeq
or
\beq{}
\frac{b_n^{(2)}}{q_{n-1}^{\t_1-1}}\geq C_1-1,
\eeq
define $a_n:=1$.
define $a_n=1$.
Because of $q_n>2^{n}$ and $\t_3<\t_2<\t_1$, there exists $n_2>n_1$ such that all these condition are not satisfied
For this $n_2$ define
\beq{}
a_{n_2}=b_{n_2}.
\eeq
So, iterating this costruction, we define $\a:=[a_1,a_2,...]$. 
By definition of $a_n$ we get that, for $n$ even:
\beq{}
\g_{-}(\a,\t_n)<\g_{+}(\a,\t_n),
\eeq
and for $n$ odd:
\beq{}
\g_{-}(\a,\t_n)>\g_{+}(\a,\t_n).
\eeq
In fact, for $n$ even we have:
\beq{}
\g(\a,\t_n)=\g_{-}(\a,\t_n)\geq C_{n-1}>\g_{+}(\a,\t_n)
\eeq
and, for $n$ odd:
\beq{}
\g(\a,\t_n)=\g_{+}(\a,\t_n)\geq C_{n-1}>\g_{-}(\a,\t_n)
\eeq
Moreover, it is easy to verify that $\t(\a)=\t$ (using remark (o)), so $\a\in D_{\bar{\t}}$ for all $\bar{\t}>\t$. 
By Lemma 2, there is a sequence $\{\bar{\t}_n\}_{n\in\Bbb{N}}$ with $\t_{n+1}<\bar{\t}_n<\t_n$ with $\a\in {\mathcal{I}}_{\bar{\t}_n}$.\qed
\nl
As an immediate consequence of Theorem B we have the following:
\cor{}
The set 
\beq{}
{\mathcal{T}}:=\Big\{\t\geq1:{\IIt}\not=\emptyset\Big\}
\eeq
is dense in $[1,+\infty).$
\ecor
\rem{}
${\IIt}=\emptyset$ for all $\t\in\Bbb{Q}$.
\erem{}
\proof
It follows by (\ref{num}). \qed
\rem{}
 ${\mathcal{I}}$ is strictly contained in $D$.
\erem{}
\proof
Define $\a:=[3,1,1,1,...]$, so $\a\in D_1$. For $\t\geq 1, n\geq 1$:
\beq{}
\su{\g_{0}(\a,\t)}=\su{\g_{0}(\a,1)}>3,
\eeq
\beq{}
\su{\gnat}=\frac{1}{q_{n}^{\t-1}}+\frac{q_{n-1}}{q_{n}^{\t}}+\frac{1}{\a_{n+2}q_{n}^{\t-1}}<\frac{3}{q_{n}^{\t-1}}
\eeq
because of $q_n<q_{n-1}$. So, for $\t\geq 1$ we have:
\beq{}
\gtal<\su{3}\leq\gtar
\eeq
Then, for all $\t\geq 1$ we have $\a\not\in{\IIt}$. \qed

\rem{d}
Given $\a\in D$, the set:
\beq{}
{\mathcal{E}}(\a):=\{\t\geq1:\a\in {\mathcal{I}}_{\t}\}
\eeq
is discrete.
\erem{}
\proof
Suppose $\t\in{\mathcal{E}}(\a)$. Let $n:=\min \{h\in\Bbb{N}_{0}:\g_h(\a,\t)=\gta\}$
Because of $g_{+}(\a,-),\g_{-}(\a,-)\in C[\t(\a),+\infty)$, it is easy to verify that there exists $\delta>0$ such that 
\beq{}
\g(\a,\t')=\g_n(\a,\t')<\g_{k}(\a,\t')
\eeq
for all $\t'\in\left(\t,\t+\delta\right)$, $k\not=n$.
If $\t=\t(\a),$ then it is clear that $\a\not\in{\mathcal{I}_{\t'}}$ for all $\t'<\t$.
If $\t>\t(\a)$, it is well defined also:
\beq{}
m:=\max\left\{h\in\Bbb{N}_{0}:\g_{h}(\a,\t)=\gta\right\}.
\eeq
Then, it is easy to check that there exists $\delta'>0$ such that:
\beq{}
\g(\a,\t')=\g_m(\a,\t')<\g_{k}(\a,\t')
\eeq
for all $\t'\in\left(\t-\delta',\t\right)$, $k\not=m$.
So, by definition of ${\IIt}$ we have $\a\not\in {\mathcal{I}_{\t'}}$ for all $\t'\in(\t-\delta',\t)\cup(\t,\t+\delta)$. \qed
\rem{}
If $\a\in D$, $\t=\t(\a)$ and there exists a strictly decreasing sequence $\{\t_n\}_{n\in\Bbb{N}}$ with $\t_n\searrow \t$ and with $\a\in{\mathcal{I}}_{\t_n}$ for all $n\in\Bbb{N}$, then $\a\not\in {\IIt}$.
\erem{}
\proof
It follows directly by Remark (\ref{d}). \qed
\nl

\section{Final observations and questions}
We have seen that, up to  an equivalent number, every Diophantine point is isolated in some Diophantine set. However, there exist Diophantine points that are always accumulation points (for example, the point defined in Remark 8). Moreover, a Diophantine number may be an isolated point for infinitely many $\t$. Indeed, by Corollary 3 it is reasonable to expect that the statement of Theorem B holds for almost every Diophantine number.
We list here some natural questions.
\begin{itemize}
    \item All the isolated points that we have construct are in $\mathcal{I}$ (i.e. the isolated point separates two intervals $I_{\g,\t}(p,q), I_{\g,\t}(a,b)$ with $p,q,a,b\in\Bbb{N}$). Is it true that for $\t>1$ the isolated points are all of this type?
    \item We have seen that $\mathcal{T}$ is dense in $[1,+\infty)$ and that $\mathcal{T}\cap\Bbb{Q}=\emptyset$.    What are the $\t\geq1$ such that ${\IIt}\not=\emptyset$? In particular, is it true that $\mathcal{T}$ is the set of Diophantine points in $[1,+\infty)$?
    \item Let $N\geq3$ and define $\D^{N}_{\g,\t}:=\{\omega\in\Bbb{R}^{N}:|\omega\cdot n|\geq\frac{\g}{|n|^{\t}}\quad\forall n\in\Bbb{Z}^{N},n\not=0\}$. What can we say about isolated points of $\D_{\g,\t}^{N}\cap\Bbb{S}^{N-1}$?
    
\end{itemize}
We have shown that in general Diophantine sets are not Cantor sets, however we believe that the following hold:
\begin{itemize}
    \item For all $\t\geq1$ there exists $\g_\t\in(0,\su{2})$ such that $\Dgt$ is a Cantor set for almost all $\g\in(0,\g_\t)$.
    \item There exists $\t^*>1$ such that, for all $\t>\t^*$:
    $$\m\left(\left\{\g\in \left(0,\su{2}\right): {\mathcal{I}}(\Dgt)\not=\emptyset\right\}\right)=0.$$
\end{itemize}
We belive also that, for any algebraic number $\a$ with degree greater then 2, there exist sequences $\t_n\searrow 1$, $\g_n\searrow 0$ such that $\a$ is an isolated point of $D_{\g_n,\t_n}$ for all $n$ (note that, if such sequences exist, by Roth Theorem $\t_n\searrow1$).

\subsubsection*{Acknowledgement}
I am very grateful to Prof. Luigi Chierchia for his suggestions, remarks and for his special support. Moreover I gratefully acknowledge useful comments by Prof. Michel Waldschmidt, Prof. Yann Bugeaud, Prof. Pappalardi, Prof. Barroero and Prof. Procesi.

\end{document}